\documentclass{article}

\usepackage{amsmath,amsthm}
\usepackage{amsrefs}

\usepackage{graphicx,psfrag}

\usepackage{subfigure}

\usepackage{wrapfig}

\theoremstyle{remark}

\theoremstyle{definition}

\theoremstyle{definition}

\begin{document}

\title{Analysis of Vibrations\\ in Large Flexible Hybrid Systems\footnote{Accepted
for publication by the \emph{Global Journal of Pure and Applied Mathematics}.
To be partially presented at the Sixth International Conference
\emph{Symmetry in Nonlinear Mathematical Physics}, June 20-26,
2005, Institute of Mathematics, National Academy of Sciences of
Ukraine, Kyiv (Kiev), Ukraine.}}

\author{Olena V. Mul\\
        \texttt{olena@mat.ua.pt}
        \and
        Delfim F.~M.~Torres\\
        \texttt{delfim@mat.ua.pt}}

\date{Department of Mathematics\\
      University of Aveiro\\
      3810-193 Aveiro, Portugal}

\maketitle

\begin{abstract}
The mathematical model of a real flexible elastic system with
distributed and discrete parameters is considered. It is a partial
differential equation with non-classical boundary conditions.
Complexity of the boundary conditions results in the impossibility
to find exact analytical solutions. To address the problem, we use
the asymptotical method of small parameter together with
the numerical method of normal fundamental systems of solutions.
These methods allow us to investigate vibrations, and a technique
for determination of complex eigenvalues of the considered boundary
value problem is developed. The conditions, at which vibration processes
of different character take place, are defined. Dependence of the vibration
frequencies on physical parameters of the hybrid system
is studied. We show that introduction of different feedbacks into
the system allow one to control the frequency spectrum, in which
excitation of vibrations is possible.
\end{abstract}

\smallskip

\noindent \textbf{Mathematics Subject Classification 2000:} 93C20, 74H10, 74H15, 74H45.

\smallskip


\smallskip

\noindent \textbf{Keywords: } hybrid dynamic control systems,
partial differential equations with non-classical boundary conditions,
vibrations, approximation of eigenvalues,
asymptotic methods, numerical methods.


\section{Introduction}

Nonlinear dynamical systems with distributed and discrete
parameters are widespread in heavy, extractive and manufacturing
industry, as well as in space-system engineering
\cite{MR80f:73002,MR88k:49009}. Different complex
dynamical processes are possible in such hybrid systems, including
vibrations, which always have negative influence on their
functioning and sometimes even can result in system breakdown.

Therefore, one of the important aims of investigation is to study
possible vibrations in such nonlinear dynamical systems with
distributed and discrete parameters, as spacecrafts of big size,
as well as to find ways how it is possible to decrease harmful
effect of vibrations on normal functioning of the systems.

Spacecrafts of big size are considered, which contain some
flexible elastic elements. Specifications for the spacecrafts (for
example, in the control of form and accuracy of space antenna
orientation) presume that intensive vibrations in such flexible
constructions are not admissible. Passive damping of vibrations is
often applied for elimination of vibrations, but strict
specifications for the spacecrafts require other methods to
control vibrations. One possibility is to use systems with active
control of vibrations. The problem is to study the effect of
active control on the dynamical characteristics of the elastic
construction with distributed parameters. There were several
attempts to solve this problem before, using mathematical models
of the system which account discrete parameters only
(see \cite{MR1835203} and references therein).
These results are not accurate enough.
Following \cite{MR1835203}, we take into account all the distributed
parameters of the flexible elastic construction. The mathematical
model becomes a nonlinear boundary value problem, which consists
of a partial differential equation and nonlinear boundary conditions.

For this new model it is possible to develop both asymptotical
and numerical methods. This allows to obtain frequencies
of possible vibrations \cite{MR1723682}.
In this way we can also determine conditions of stability
under active control of vibrations.
After analysis of influence of different parameters of the
considered systems it is possible to make a conclusion about the
optimal parameters.

The results of this paper may be applied in space-system
engineering for design of new improved technical systems.


\section{Mathematical Model}

In the given paper a flexible elastic system
with distributed and discrete parameters is analyzed.

\subsection{Basic Mathematical Model}

The basic mathematical model, which is used for investigations,
is borrowed from \cite{MR1835203}:

\begin{equation}
\label{eq:4} \rho\,S  \frac{\partial^2 u}{\partial t^2} -ES
\frac{\partial^2 u}{\partial x^2} =ES\beta\, \frac{\partial^3
u}{\partial x^2 \partial t} \, ,
 \end{equation}

\begin{equation}
\label{eq:5} x = 0 : \quad u \left( 0,t \right) =0 \, ,
 \end{equation}

\begin{multline}
\label{eq:6}
x = l : \quad
mES\beta\,  \frac{\partial^4 u}{\partial x \partial t^3}
+ m(d+b) \frac{\partial^3 u}{\partial t^3}
 + E S(m + \beta b) \frac{\partial^3 u}{\partial x \partial t^2}\\
+mc \frac{\partial^2 u}{\partial t^2} + E S(b + \beta\,c)
\frac{\partial^2 u}{\partial x \partial t} +ESc \frac{\partial
u}{\partial x}  =0 \, ,
 \end{multline}
where $u(x,t)$ is the displacement of a point of the bar with abscissa
$x$, $S$ the bar cross-section area, $\rho$ the unit volume mass,
$E$ the modulus of the material elasticity, $\beta$ the coefficient
of dissipation in the bar material, $b$ the damping factor of
the executive mechanism, $c$ the rigidity of the centering spring,
$d$ the feedback coefficient, $m$ the specimen mass,
and $l$ is the length of the elastic bar.

It is necessary to point out that problem \eqref{eq:4}--\eqref{eq:6}
is a nonconservative boundary value problem. This fact introduces
considerable difficulties for the computation of an exact analytical
solution of the problem. The non-conservatism is caused by the presence
of the odd derivatives on time both in the equation of motion \eqref{eq:4}
and in the boundary conditions \eqref{eq:5} and \eqref{eq:6}.
This allows to assume that eigenvalues of the considered boundary
value problem are complex numbers, what is in agreement with the results
obtained in \cite{MR88k:49009,MR2081767,MR1383822}.

Problem \eqref{eq:4}--\eqref{eq:6} is a nonclassical
boundary value problem, and it is necessary to take into
account real features of objects functioning.
The aim is to obtain approximations for the eigenvalues
of this system. We will tackle the problem both
from a theoretical and a numerical point of view.

\subsection{Dimensionless Form of the Mathematical Model }

In order to account the influence of different system parameters
onto possible vibrations, as well as to simplify cumbersome
calculations, we introduce dimensionless variables
$\bar{u} = \frac{u}{l}$,
$\bar{x} = \frac{x}{l}$,
$\tau = \frac{ta}{l}$, and dimensionless parameters
\begin{equation}
\label{eq:8} \varepsilon_1 = \frac{\beta a}{l} \, , \quad \mu =
\frac{b a}{E S} \, , \quad \nu = \frac{d a}{E S} \, , \quad \eta =
\frac{m}{\rho S l} \, , \quad \delta = \frac{E S}{c l} \, ,
\end{equation}
where $a^2 = \frac{E}{\rho}$.
Then the dimensionless form of boundary problem
\eqref{eq:4}--\eqref{eq:6} may be written as follows:

\begin{equation}
\label{eq:9}
\frac{\partial^2 \bar{u}}{\partial \tau^2}
 - \frac{\partial^2 \bar{u}}{\partial \bar{x}^2}
={\it \varepsilon_1}\, \frac{\partial^3 \bar{u}}{\partial
\bar{x}^2 \partial \tau} \, ,
\end{equation}

\begin{equation}
\label{eq:10} \bar{x} = 0 : \quad {\it \bar{u}} \left( 0,\tau
\right) =0 \, ,
\end{equation}

\begin{multline}
\label{eq:11} \bar{x} = 1 : \quad {\it
\varepsilon_1}\,\eta\,\delta\, \frac{\partial^4 \bar{u}}{\partial
\bar{x} \partial \tau^3} +\eta\,\delta\, (\nu + \mu)
\frac{\partial^3 \bar{u}}{\partial \tau^3}
 + \delta (\eta + \varepsilon_1 \,\mu)
\frac{\partial^3 \bar{u}}{\partial \bar{x} \partial \tau^2} \\
+\eta\, \frac{\partial^2 \bar{u}}{\partial \tau^2}
 +(\mu\,\delta + \varepsilon_1)
\frac{\partial^2 \bar{u}}{\partial \bar{x} \partial \tau} +
\frac{\partial \bar{u}}{\partial \bar{x}} =0 \, .
\end{multline}

Mathematical model \eqref{eq:9}--\eqref{eq:11} describes
dynamical processes in the system with distributed parameters,
and such system has infinite number of degrees of freedom.
Therefore it is natural to assume that stationary operating
conditions in the considered system are not always single-frequency ones.


\section{Asymptotical Approach}

The first approach, which may be used for problem solving, is an
asymptotical one.
Although at present numerous efficient methods for the investigation of vibrations
have been developed, and mathematically justified, the most efficient ones are the asymptotical
methods, especially the averaging method \cite{MR2091174}. This method is often used
for the solution of boundary value problems appearing in the simulation of the behavior
of real processes and in problems of optimal control.

The idea of the averaging method is as follows: using some special operator, one replaces the system
of differential equations under investigation by another system, the so-called averaged system.
The last one should be simpler than the original system but at the same time it should describe
the main features of the phenomenon under investigation.

Today the averaging method is developed for various classes of differential equations
with large and small parameters, although in most cases the problem of its justification
remained unsolved till now.

In case of the method of the small parameter, it is necessary to choose systems of equations
of higher approximations whose solutions approximate the solutions of the original system
of equations to within values proportional to integer powers of some small parameter.

\subsection{Basic Assumption}

In order to solve the boundary value problem \eqref{eq:9}--\eqref{eq:11},
it is expedient to apply the modifications of the asymptotical method
of the small parameter.

To do this, it is necessary to take into account
that equations \eqref{eq:9}--\eqref{eq:11}
contain five dimensionless parameters \eqref{eq:8}.
These parameters, as a rule, are less than unit,
although they have different orders of smallness.
Dimensionless parameters $\varepsilon_1$ and $\mu$ are
connected, respectively, with the viscous friction in the material
of the bar and the centering spring, whereas dimensionless
parameter $\nu$ is connected with the coefficient
of the velocity feedback. Therefore, parameters
$\varepsilon_1$, $\mu$ and $\nu$ are, for real dynamical systems,
of magnitude order less than dimensionless parameters $\eta$ and $\delta$.
This allow us to assume that dimensionless parameters characterizing
dissipation of different physical nature are proportional to some
common small parameter $\varepsilon$, \textrm{i.e.}

\begin{equation*}
\label{eq:12} \varepsilon_1 = \varepsilon \, \bar{\varepsilon}_1
\, , \quad \mu = \varepsilon \, \bar{\mu} \, , \quad \nu =
\varepsilon \, \bar{\nu} \, .
\end{equation*}
With such an assumption, mathematical model
\eqref{eq:9}--\eqref{eq:11} may be transformed to the form:

\begin{equation}
\label{eq:9-nc} \frac{\partial^2 \bar{u}}{\partial \tau^2}
 - \frac{\partial^2 \bar{u}}{\partial \bar{x}^2}
={\it \varepsilon \, \bar{\varepsilon}_1}\, \frac{\partial^3
\bar{u}}{\partial \bar{x}^2 \partial \tau} \, ,
\end{equation}

\begin{equation}
\label{eq:10-nc} \bar{x} = 0 : \quad {\it \bar{u}} \left( 0,\tau
\right) =0 \, ,
\end{equation}

\begin{multline}
\label{eq:11-nc} \bar{x} = 1 : \quad {\it \varepsilon \,
\bar{\varepsilon}_1}\,\eta\,\delta\, \frac{\partial^4
\bar{u}}{\partial \bar{x} \partial \tau^3} + \varepsilon \,
\eta\,\delta\, ( \bar{\nu} + \bar{\mu}) \frac{\partial^3
\bar{u}}{\partial \tau^3}
 + \delta (\eta + \varepsilon^{2} \, \bar{\varepsilon}_1 \, \bar{\mu})
\frac{\partial^3 \bar{u}}{\partial \bar{x} \partial \tau^2} \\
+\eta\, \frac{\partial^2 \bar{u}}{\partial \tau^2}
 +\varepsilon \, (\bar{\mu}\,\delta + \bar{\varepsilon}_1)
\frac{\partial^2 \bar{u}}{\partial \bar{x} \partial \tau} +
\frac{\partial \bar{u}}{\partial \bar{x}} =0 \, .
\end{multline}
We are now in conditions to apply for investigations some modifications
of asymptotical methods of the small parameter.


\subsection{Conservative Case}

As it was mentioned above, mathematical model \eqref{eq:9-nc}--\eqref{eq:11-nc}
includes different factors, which characterize energy dissipation in the system.
Therefore both equation of motion \eqref{eq:9-nc} and boundary conditions
\eqref{eq:10-nc}, \eqref{eq:11-nc} contain the odd derivatives on time, which
represent non-conservatism of the physical model of our hybrid system.

But for the first approximation we can neglect all the dissipative
factors of the mathematical model, and consider its conservative case,
\textrm{i.e.} put $\varepsilon = 0$.
Then we will have $\varepsilon_1 = \mu = \nu = 0$ and problem
\eqref{eq:9-nc}--\eqref{eq:11-nc} simplifies to

\begin{equation}
\label{eq:9:ep0}
\frac{\partial^2 \bar{u}}{\partial \tau^2}
 - \frac{\partial^2 \bar{u}}{\partial \bar{x}^2} =0 \, ,
\end{equation}

\begin{equation}
\label{eq:10:ep0}
\bar{x} = 0 : \quad
{\it \bar{u}} \left( 0,\tau \right) =0 \, ,
\end{equation}

\begin{equation}
\label{eq:11:ep0} \bar{x} = 1 : \quad \delta \eta \frac{\partial^3
\bar{u}}{\partial \bar{x} \partial \tau^2} +\eta\,
\frac{\partial^2 \bar{u}}{\partial \tau^2} + \frac{\partial
\bar{u}}{\partial \bar{x}} =0 \, .
\end{equation}

It is easy to solve exactly problem \eqref{eq:9:ep0}--\eqref{eq:11:ep0} and
obtain the next transcendental equation for determination of the eigenvalues:

\begin{equation}
\label{eq:17epsilon0}
\cot \left( \omega \right) = \frac{\eta\,\omega}{1-\eta\,\delta\,\omega^{2}} \, ,
\end{equation}
where $\omega$ is an eigenvalue of the conservative boundary
problem, which is a real number.

Formula \eqref{eq:17epsilon0} can be used as the initial estimate of the
complex eigenvalues of the original nonconservative boundary problem.


\subsection{First Asymptotical Method}

According with \cite{myNY98}, we will look for the solution
of equation \eqref{eq:9-nc} in the following form:

\begin{equation}
\label{eq:16}
\bar{u}(\bar{x},\tau) =
\left(A\cos\left(\lambda\,\bar{x}\right)
+B\sin\left(\lambda\,\bar{x}\right)\right)
\textrm{e}^{\left(-\frac{1}{2}\,\varepsilon\,\bar{\varepsilon}_1\,\lambda^{2}
+i\lambda \right) \tau} \, ,
\end{equation}
where $\lambda$ is an eigenvalue of the nonconservative boundary
problem, $A$ and $B$ are arbitrary constants, and $i=\sqrt{-1}$ .

Solution \eqref{eq:16} satisfies equation \eqref{eq:9-nc} with
accuracy $\varepsilon^{2}$ . After substitution of solution
\eqref{eq:16} into boundary conditions \eqref{eq:10-nc} and
\eqref{eq:11-nc}, and their satisfaction with the same accuracy
$\varepsilon^{2}$, we obtain a transcendental equation of the form

\begin{equation}
\label{eq:17} \cot \left( \lambda \right) \approx
\frac{\eta\,\lambda\, + i \, \varepsilon \, \eta\,\lambda^2\
\left( {\it \bar{\varepsilon}_1}+{\it \bar{\mu}}\,\delta +{\it
\bar{\nu}}\,\delta \right)}{1-\eta\,\delta\,\lambda^{2} + i
\varepsilon \lambda \left(\bar{\varepsilon}_1
-2\eta\,\delta\,\lambda^{2}\,\bar{\varepsilon}_1
+\bar{\mu}\,\delta\right)} \, .
\end{equation}
The solution of this equation at $\varepsilon=0$ is known, and it determines
eigenvalues $\omega$ of the corresponding conservative boundary problem
\eqref{eq:9:ep0}--\eqref{eq:11:ep0}.

In solving transcendental equation \eqref{eq:17}, it is necessary to account
smallness of damping forces. This allow us to assume that the correction factor
for non-conservatism is a quantity of order $\varepsilon$, \textrm{i.e.}

\begin{equation}
\label{eq:18} \lambda = w + \varepsilon \lambda_1 \, ,
\end{equation}
where $w$ is an eigenvalue of the corresponding conservative problem,
and $\varepsilon \lambda_1$ is the correction factor for the
non-conservatism of the problem.

Substituting \eqref{eq:18} into transcendental equation \eqref{eq:17},
and cutting off a Taylor series for left-hand side of this equation as follows

\begin{equation*}
\cot\left(w+\varepsilon \lambda_1\right) = \cot \left( w \right) -
\left( 1+  \cot \left( w \right)^{2} \right) {\it
\lambda_1}\,\varepsilon+O \left( {\varepsilon}^{2}\right) \, ,
\end{equation*}
it is possible to present the next expression for $\lambda_1$:

\begin{equation*}
\label{eq:19} {\it \lambda_1} \approx {\frac
{i\,\eta\,\delta\,{w}^{2}\left[ \left(\delta \,{\it
\bar{\mu}}+{\it \bar{\nu}}\,\delta-{\it
\bar{\varepsilon}}_1\right) \,\eta\,{w}^{2} -{\it \bar{\nu}}
\right]}{{\eta}^{2}{\delta}^{2}{w}^{4}+
\left(\delta\,{\eta}-2\,\delta+{\eta}\right)\eta {w}^{2} +
\eta+1}} \, .
\end{equation*}

Thus, eigenvalues of the nonconservative problem
\eqref{eq:9-nc}--\eqref{eq:11-nc}
are approximately determined in the following way:

\begin{equation}
\label{eq:20}
\lambda \approx w + {\frac {i\, \varepsilon \,
\eta\,\delta\,{w}^{2}\left[ \left(\delta \,{\it \bar{\mu}}+{\it
\bar{\nu}}\,\delta-{\it \bar{\varepsilon}_1}\right)
\,\eta\,{w}^{2} -{\it \bar{\nu}}
\right]}{{\eta}^{2}{\delta}^{2}{w}^{4}+
\left(\delta\,{\eta}-2\,\delta+{\eta}\right)\eta {w}^{2} +
\eta+1}} \, .
\end{equation}

Now, it is possible to substitute the complex value found for $\lambda$ into
the exponential function of solution \eqref{eq:16}.
Analysis of the real part of the exponent allows to obtain condition

\begin{equation}
\label{eq:23}
{\frac {{\it \varepsilon_1}{\eta}^{2}{\delta}^{2}{w}^{4}
+  \left( 2\,{\eta}^{2}{\delta}^{2}(\mu + \nu)+{\it
\varepsilon_1}{\eta}^{2}(1-\delta) -2\,\eta\,\delta\,{\it
\varepsilon_1} \right) {w}^{2}+{\it \varepsilon_1} (1 +
\eta)-2\,\eta\,\delta\,\nu}{{\eta}^{2}{\delta}^{2}{w}^{4 }+
\left[\eta(1 + \delta) -2\,\delta \right]\eta {w}^{2}+\eta+1}} \leq 0 \, ,
\end{equation}
whose satisfaction leads to self-excitation of vibrations, and
self-vibrations in case of equality.

It is necessary also to point out that the imaginary part of the
exponent represents the fundamental frequency of the
nonconservative system.

Parameter $\nu$, which may be adjusted by introducing different
feedbacks, enters into inequality \eqref{eq:23}. Hence, the
possibility appears to control the frequency spectrum, on which
vibrations excitation is possible, by varying dissipative
characteristics of the hybrid dynamical system.


\subsection{Second Asymptotical Method}

It is known that in a conservative system (at $\varepsilon=0$) periodic motions appear
with any amplitude dependent on initial conditions. But in the nonconservative system
(at $\varepsilon\neq0$) they have only fixed amplitudes, which conform to equality of an
energy present at the expense of negative resistance and energy dissipation, which is always
present in a real dynamical system.

Thus, a solution of nonconservative boundary problem
\eqref{eq:9-nc}--\eqref{eq:11-nc} may be sought for in the form

\begin{equation}
\label{eq:56} \bar{u}(\bar{x},\tau) =
\bar{u}_0\left(\bar{x},\tau\right)+
\varepsilon\,\bar{u}_1\left(\bar{x},\tau\right) \, ,
\end{equation}
where $\bar{u}_0\left(\bar{x},\tau\right)$ is one of the solutions
of the conservative problem \eqref{eq:9:ep0}--\eqref{eq:11:ep0}.

It is easy to obtain the next expression for such solution:

\begin{equation}
\label{eq:57}
\bar{u}_0\left(\bar{x},\tau\right) = A \sin(w \bar{x})\cos(w\tau) \, ,
\end{equation}
where $A$ has the meaning of amplitude, and it is necessary to select
it close to amplitude of boundary cycle of self-sustained
vibration system.

Then, after taking into account expressions \eqref{eq:56} and
\eqref{eq:57}, as well as relationship \eqref{eq:17epsilon0}
for eigenvalues $\omega$ of conservative boundary problem
\eqref{eq:9:ep0}--\eqref{eq:11:ep0}, the boundary problem
\eqref{eq:9-nc}--\eqref{eq:11-nc} may be written as follows:

\begin{equation}
\label{eq:510} \frac{\partial^2 \bar{u}_1}{\partial \tau^2} -
\frac{\partial^2 \bar{u}_1}{\partial \bar{x}^2} =
\bar{\varepsilon}_1 A w^3 \sin(w \bar{x}) \sin(w\tau) \, ,
\end{equation}

\begin{equation}
\label{eq:511} \bar{x} = 0 : \quad \bar{u}_1(0,\tau) = 0 \, ,
\end{equation}
\begin{multline}
\label{eq:512} \bar{x} = 1 : \  A \omega \cos(\omega\tau) \left[
(\eta\delta \omega^2-1) \cos(\omega)+ \eta \omega \sin(\omega) \right]
 = \varepsilon \biggl[ \eta \delta \frac{\partial^3
\bar{u}_1}{\partial \bar{x} \partial \tau^2} + \eta
\frac{\partial^2 \bar{u}_1}{\partial \tau^2}
\\ + \frac{\partial \bar{u}_1}{\partial \bar{x}}
+ A \omega^2 \sin(\omega\tau) \left[ \left( \bar{\varepsilon}_1 (\eta \delta
\omega^2 -1) - \bar{\mu} \delta \right) \cos(\omega) + \eta \delta \omega
(\bar{\nu}+ \bar{\mu}) \sin(\omega) \right] \biggr] \, .
\end{multline}

Boundary problem \eqref{eq:510}--\eqref{eq:512} may be considered
as a mathematical model of the conservative system with natural frequency $\omega$,
on which an external force of resonance frequency effects. It is known \cite{myDresd}
that a stationary solution exists only on the condition of orthogonality of external
force and normal mode of a system. This condition determines an amplitude,
which does not increase with time.

Therefore, for the considered boundary problem
it is expedient to choose the solution in the form
\begin{equation*}
\label{eq:515}
\bar{u}_1(\bar{x},\tau) = U(\bar{x}) \sin(\omega\tau) \, ,
\end{equation*}
which reduces our problem to ordinary differential equations:

\begin{equation}
\label{eq:516} \frac{d^2 U}{d\bar{x}^2} + \omega^2 U = -
\bar{\varepsilon}_1 A \omega^3 \sin(\omega\bar{x}) \, ,
\end{equation}
\begin{equation}
\label{eq:517}
\bar{x} = 0 : \quad U = 0 \, ,
\end{equation}
\begin{multline}
\label{eq:518} \bar{x} = 1 : \quad A \omega \cos(\omega\tau) \left[
(\eta\delta \omega^2 -1) \cos(\omega) + \eta \omega \sin(\omega) \right] \\
= \varepsilon \sin(\omega\tau) \biggl[A \omega^3 \eta \delta (\bar{\nu} +\bar{\mu})
\sin(\omega) + A \omega^2 (\bar{\varepsilon}_1 (\eta\delta \omega^2 - 1)
-\bar{\mu}\delta )\cos(\omega) \\
- \frac{dU}{d\bar{x}} (\eta\delta \omega^2 - 1) - \eta U \omega^2\biggr] \, .
\end{multline}

Searching the general solution of equation \eqref{eq:516} in the form

\begin{equation*}
\label{eq:29p}
U (\bar{x}) = (B_{1}+B_{2}\bar{x})\cos(\omega\bar{x}) +
              (C_{1}+C_{2}\bar{x})\sin(\omega\bar{x}) \, ,
\end{equation*}
where $B_{1}$, $B_{2}$, $C_{1}$, $C_{2}$ are arbitrary constants, and satisfying
boundary conditions \eqref{eq:517} and \eqref{eq:518} allows to obtain the next condition:

\begin{multline}
\label{eq:523} {\frac {{\it
C_2}\left(\left[{\eta}^{2}{\omega}^{2}(1+\delta)
+(\eta\delta{\omega}^{2}-1)^{2}-\eta\right]\omega \cot \left(\omega{\it \bar{x}}
\right) - (\eta\delta{\omega}^{2}-1)^{2}\right)}{ \left[ {\it
\varepsilon_1}({\eta}^{2}{\omega}^{2}(1-\delta)+(\eta\delta{\omega}^{2}-1)^{2}+\eta)
+2({\eta}^{2}{\omega}^{2}{\delta}^{2}({\it \bar{\nu}}+{\it
\bar{\mu}})-\eta{\it \bar{\nu}}\delta )\right]{\omega}^{3}}} \leq 0 \, .
\end{multline}

Inequality \eqref{eq:523} gives a condition for self-vibration excitation
in the considered large flexible hybrid system. In the case of equality
in \eqref{eq:523} we have the most interesting case of
self-excited vibrations in the system.

The presence in \eqref{eq:523}
of arbitrary constant $C_{2}$ signifies the dependence of vibrations amplitude $A$
on the initial conditions, which were not considered in our boundary problem.

Thus, it was confirmed again \cite{MR88k:49009} that, in order to find some
possibilities of vibrations control in the hybrid system,
one does not need to introduce the initial conditions,
being enough to study the frequencies of vibrations.


\section{Numerical Approach}

The other approach, which may be used for problem investigations,
is a numerical one.

In this paper, the numerical method of the normal
fundamental system of solutions is used for determination of the
eigenvalues of the boundary value problem \eqref{eq:9}--\eqref{eq:11}.
The justification of the convergence of the method of normal fundamental
systems of solutions, together with the theory which explains its good
accuracy when applied to vibration problems, can be found
in the monograph \cite{Norm-fund}.

\subsection{Basic Method}

Taking into account that the eigenvalues are complex numbers,
the solution may be searched in the following form:

\begin{equation}
\label{eq:nm:formSol}
u\left(\bar{x},\tau\right) = \left( u_1 \left( \bar{x} \right) +i u_2 \left( \bar{x} \right)  \right)
 {\mathrm{e}^{ \left( q+i\omega \right) \tau}} \, ,
\end{equation}
where $q$ and $\omega$ are, respectively, the real and the imaginary parts
of the eigenvalues. They are real constants, which should be determined.

After substitution of \eqref{eq:nm:formSol} in the boundary problem
\eqref{eq:9}--\eqref{eq:11}, the real and imaginary parts of the equation
of motion \eqref{eq:9} are given by
\begin{equation}
\label{eq:nm:mot-ri}
\begin{split}
\left( {q}^{2}- {\omega}^{2}\right)u_1
-2\, q \omega u_2-
\left(1 + \varepsilon_1 q \right)
 u''_1 + \varepsilon_1\, \omega u''_2 & = 0 \, , \\
2 q \omega u_1 + \left(q^2 - \omega^2\right) u_2
- \varepsilon_1 \omega u''_1
- \left(\varepsilon_1 q + 1\right) u''_2 & = 0 \, .
\end{split}
\end{equation}

In order to ensure the compactness of further expressions,
let us make the following definitions:
\begin{equation*}
\label{eq:const}
\begin{split}
D_1 & \doteq \eta (q^2 - \omega^2) + \eta \delta q^3 (\nu + \mu) - 3 \eta \delta q \omega^2 (\nu + \mu) \, ,\\
D_2 & \doteq -3 \eta \delta q^2 \omega (\nu + \mu) + \eta \omega \left(\mu \delta \omega^2 - 2q + \nu \delta \omega^2\right) \, ,\\
D_3 & \doteq \delta (q^2 - \omega^2) \left(\varepsilon_1 \mu + \eta\right)
           + \varepsilon_1 \eta\delta q \left(q^2-3\omega^2\right) + q (\mu\delta + \varepsilon_1)+1 \, , \\
D_4 & \doteq - \omega (\varepsilon_1 + \mu\delta)+\varepsilon_1\eta\delta \omega (\omega^2 - 3 q^2)
           -2\delta q \omega (\eta + \varepsilon_1 \mu) \, .
\end{split}
\end{equation*}
With these notations, after separation of their real and imaginary parts,
boundary conditions \eqref{eq:10} and \eqref{eq:11} can be written as

\begin{equation}
\label{eq:bc1}
\begin{split}
\bar{x} = 0 : \quad u_1 & = 0 \, , \\
u_2 & = 0 \, ;
\end{split}
\end{equation}
\begin{equation}
\label{eq:bc2}
\begin{split}
\bar{x} = 1 : \quad D_1 u_1 + D_2 u_2 + D_3 u'_1 + D_4 u'_2 & = 0 \, , \\
- D_2 u_1 + D_1 u_2 - D_4 u'_1 + D_3 u'_2 & = 0 \, .
\end{split}
\end{equation}

An application of the numerical method of normal fundamental systems of
solutions demands a presentation of the problem as a system of ordinary differential
equations of the first order in the normal form, which satisfy some boundary conditions.

Therefore, in order to reduce \eqref{eq:nm:mot-ri} to the necessary form,
let us enter new functions
\begin{equation}
\label{eq:nm:new-funct}
\gamma_1 = u_1 \, , \quad \gamma_2 = u_2 \, , \quad
\gamma_3 = u'_1 \, , \quad \gamma_4 = u'_2 \, .
\end{equation}

Then, the desired system of differential equations in the normal form may be written as
\begin{equation}
\label{eq:nm:norm-syst}
\begin{cases}
\gamma'_1 = \gamma_3 \, , \\
\gamma'_2 = \gamma_4 \, , \\
\gamma'_3 = K_1 \gamma_1 - K_2 \gamma_2 \, , \\
\gamma'_4 = K_2 \gamma_1 + K_1 \gamma_2 \, ,
\end{cases}
\end{equation}
where coefficients $K_1$ and $K_2$ are dependent on the sought-for quantities $q$ and $\omega$:
\begin{equation*}
K_1 = \frac{q^2-\omega^2+\varepsilon_1 q \left(q^2 + \omega^2\right)}{{\varepsilon_1}^{2}{
\omega}^{2}+\left(1+\varepsilon_1\,q\right)^2} \, , \quad
K_2 = \frac{\left(2q+\varepsilon_1\left(q^2 + \omega^2\right)\right)\omega}{{\varepsilon_1}^{2}{
\omega}^{2}+\left(1+\varepsilon_1\,q\right)^2} \, .
\end{equation*}

In new functions \eqref{eq:nm:new-funct}, the boundary conditions
\eqref{eq:bc1} and \eqref{eq:bc2} will be transformed to the next ones:

\begin{equation}
\label{eq:bc1:gamma}
\bar{x} = 0 : \quad \gamma_1 = \gamma_2  = \gamma_3 = \gamma_4 = 0 \, ,
\end{equation}
\begin{equation}
\label{eq:bc2:gamma}
\begin{split}
\bar{x} = 1 : \quad
D_1 \gamma_1 + D_2 \gamma_2 + D_3 \gamma_3 + D_4 \gamma_4 = 0 \, ,\\
-D_2 \gamma_1 + D_1 \gamma_2 - D_4 \gamma_3 + D_3 \gamma_4 = 0 \, .
\end{split}
\end{equation}

Now it is necessary to solve the boundary value problem \eqref{eq:nm:norm-syst},
\eqref{eq:bc1:gamma}, \eqref{eq:bc2:gamma}. By means of any known numerical method,
one can solve the Cauchy problem for the set of equations \eqref{eq:nm:norm-syst}
four times with the following initial conditions:

\begin{equation*}
\gamma_{j,k}(0,\omega,q) =
\begin{cases}
1 & \text{ if } k = j \\
0 & \text{ if } k \ne j
\end{cases} \, , \quad j,\, k = 1,2,3,4 \, ,
\end{equation*}
where the first subscript of the function $\gamma$ is a solution number, and the second
one is a function number.

Thus, the normal fundamental system of solutions
$\left\{\gamma_{j,k}(\bar{x},\omega,q)\right\}_{j,k=1,2,3,4}$ for the set of differential
equations \eqref{eq:nm:norm-syst} may be generated. With the help of such system of solutions,
general solution $\gamma_k(\bar{x},\omega,q)$
of set \eqref{eq:nm:norm-syst} may be written as follows:

\begin{equation*}
\gamma_k(\bar{x},\omega,q) = \sum_{j=1}^{4} C_{j} \gamma_{j,k}\left(\bar{x},\omega,q\right) \, ,
\end{equation*}
where $C_j$ are some coefficients, which should be determined.

Using the property of the normal fundamental system of solutions
\begin{equation*}
\label{eq:property}
C_k =  \gamma_k(0,\omega,q), \quad k = 1,2,3,4 \, ,
\end{equation*}
from boundary condition \eqref{eq:bc1:gamma} it is easy to conclude that $C_1 = C_2 = 0$.

Now, satisfying also boundary conditions \eqref{eq:bc2:gamma}, we will obtain
two linear algebraic equations in coefficients $C_3$ and $C_4$:
\begin{equation}
\label{eq:lin:alg:set}
\begin{cases}
C_3 E_1(1,\omega,q) + C_4 E_2(1,\omega,q) = 0 \, ,\\
C_3 E_3(1,\omega,q) +C_4 E_4(1,\omega,q) = 0 \, ,
\end{cases}
\end{equation}
where coefficients $E_s(1,\omega,q)$ $(s = 1,2,3,4)$ are defined by
\begin{equation*}
\begin{split}
E_1(1,\omega,q) &\doteq D_1\gamma_{3,1}(1) + D_2\gamma_{3,2}(1)
+D_3\gamma_{3,3}(1)+D_4\gamma_{3,4}(1) \, ,\\
E_2(1,\omega,q) &\doteq D_1\gamma_{4,1}(1)+D_2\gamma_{4,2}(1)
+D_3\gamma_{4,3}(1)+D_4\gamma_{4,4}(1)\, ,\\
E_3(1,\omega,q) &\doteq -D_2\gamma_{3,1}(1) + D_1\gamma_{3,2}(1)-D_4\gamma_{3,3}(1)+D_3\gamma_{3,4}(1)\, ,\\
E_4(1,\omega,q) &\doteq -D_2\gamma_{4,1}(1)+D_1\gamma_{4,2}(1)-D_4\gamma_{4,3}(1)+D_3\gamma_{4,4}(1) \, .
\end{split}
\end{equation*}

Then, the relation

\begin{equation}
\label{eq:determ}
\Delta(\omega,q) =
\left|
\begin{array}{cc}
  E_1(1,\omega,q) & E_2(1,\omega,q) \\
  E_3(1,\omega,q) & E_4(1,\omega,q) \\
\end{array}
\right| = 0
\end{equation}
is a necessary and sufficient condition for the existence of a non-trivial solution
of the homogeneous set of linear algebraic equations \eqref{eq:lin:alg:set}.

This relation is also a necessary and sufficient condition for the existence of a
non-trivial solution of the boundary value problem \eqref{eq:4}, \eqref{eq:5}, \eqref{eq:6}.
It is an equation for a pair of numbers $q$ and $\omega$, which may be determined
as its roots. Thus, from equation \eqref{eq:determ} it is possible to determine complex
eigenvalues $q+ i \omega$ of the boundary value problem.

It is necessary also to point out that the determinant $\Delta(\omega,q)$ from expression
\eqref{eq:determ} is always non-negative, which is in agreement with the results obtained
in \cite{MR88k:49009}.

\subsection{Improved Method}

It is obvious that the accuracy of determination of the boundary problem eigenvalues
depends on the length of the integration interval for system \eqref{eq:nm:norm-syst}.
Therefore, in order to decrease the calculating error, it is necessary to divide
the integration interval $[0,1]$ into some smaller intervals by points $\bar{x}=\bar{x}_{i}$,
where $i=\overline{1,n-1}$, assuming $\bar{x}_{0}=0$, $\bar{x}_{n}=1$.

Then solutions $\gamma_k(\bar{x},\omega,q)$ of boundary value problem \eqref{eq:nm:norm-syst},
\eqref{eq:bc1:gamma}, \eqref{eq:bc2:gamma} may be formed on corresponding subintervals
by solutions $\gamma_k^{(i)}(\bar{x},\omega,q)$ of systems of equations

\begin{equation}
\label{eq:nm:norm-syst-div}
\begin{cases}
[\gamma ^{(i)}_1]' = \gamma_3^{(i)} \, , \\
[\gamma ^{(i)}_2]' = \gamma_4^{(i)} \, , \\
[\gamma ^{(i)}_3]'= K_1 \gamma_1^{(i)} - K_2 \gamma_2^{(i)} \, , \\
[\gamma ^{(i)}_4]'= K_2 \gamma_1^{(i)} + K_1 \gamma_2^{(i)} \, ,
\end{cases}
\end{equation}
where function $\gamma$ subscript is a solution number, and its superscript is a
subinterval number, $i=\overline{1,n}$.

The solutions $\gamma_k^{(i)}(\bar{x},\omega,q)$ of \eqref{eq:nm:norm-syst-div}
must satisfy boundary conditions \eqref{eq:bc1:gamma} and \eqref{eq:bc2:gamma}
and, besides, some conjugation conditions, which link them at the ends
of the subintervals. Having in mind that points
$\bar{x}_{i}$ of a partition of the integration interval $[0,1]$ are arbitrary,
we can write the following conjugation conditions:

\begin{equation*}
\label{eq:conjug:cond}
\begin{split}
\gamma_k^{(i)}(\bar{x}_{i},\omega,q) = \gamma_k^{(i+1)}(\bar{x}_{i},\omega,q) , \\
[\gamma_k^{(i)}(\bar{x}_{i},\omega,q)]' = [\gamma_k^{(i+1)}(\bar{x}_{i},\omega,q)]' ,
\end{split}
\end{equation*}
where $i=\overline{1,n-1}$. In this way, solutions $\gamma_k(\bar{x},\omega,q)$
of boundary value problem \eqref{eq:nm:norm-syst}, \eqref{eq:bc1:gamma},
\eqref{eq:bc2:gamma} are defined as follows:

\begin{equation}
\label{eq:solut:int}
\gamma_k(\bar{x},\omega,q) = \gamma_k^{(i)}(\bar{x},\omega,q) ,
\end{equation}
where $\bar{x}\in (\bar{x}_{i-1},\bar{x}_{i}]$, $i=\overline{1,n}$.

So now it is possible to use the above described technique of solution
for all the systems of equations \eqref{eq:nm:norm-syst-div}. Then, accounting
\eqref{eq:solut:int}, we get a necessary and sufficient condition for the existence
of the non-trivial solution of the considered problem:

\begin{equation}
\label{eq:determ:div}
\Delta(\omega,q) =
\left|
\begin{array}{cc}
  E_1(\bar{x}_{n},\omega,q) & E_2(\bar{x}_{n},\omega,q) \\
  E_3(\bar{x}_{n},\omega,q) & E_4(\bar{x}_{n},\omega,q) \\
\end{array}
\right| = 0 \, ,
\end{equation}
where $E_1$, $E_2$, $E_3$, $E_4$ can be expressed by recurrence formulas.

Again, function $\Delta(\omega,q)$ from equation \eqref{eq:determ:div} is always non-negative.
This fact complicates the numerical determination of the eigenvalues and, in this connection,
we have used the method of minimum direct search for determination of the pairs $q$ and $\omega$,
which minimize function $\Delta(\omega,q)$.


\section{Discussion of the Results and Conclusions}

In this paper we have used both asymptotical and numerical methods,
successfully solving nonconservative boundary value problems which
describe a class of hybrid control dynamical systems appearing
in applications of space-engineering design.

We present now some results of numerical calculations for
the boundary value problem \eqref{eq:9}--\eqref{eq:11} as well as their analysis.

First, let us consider the results of calculations by the first asymptotical method.
According to \eqref{eq:17epsilon0} and \eqref{eq:20}, the first two fundamental
frequencies of the nonconservative system are the next ones:
\begin{equation}
\label{eq:w}
\omega_1 = 0.3534042288, \quad  \omega_2 = 2.904816694 \, .
\end{equation}

\begin{figure}
  \centering
   \psfrag{nu}{$\nu$}
   \psfrag{w}{\hspace*{-0.3cm}$\omega$}
   \includegraphics[scale=0.4]{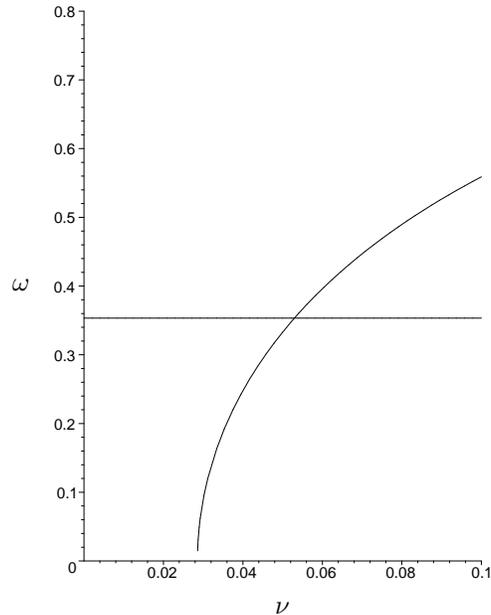}
\caption{Asymptotical Methods: $\varepsilon_1 = 0.005$, $\mu = 0.008$, $\eta=7$, $\delta=0.1$}
\label{fig:AM1:w-nu}
\end{figure}

Now, analysis of condition \eqref{eq:23} of vibrations excitation shows that
a boundary value of frequency always exists in the considered hybrid system,
which separates frequencies domains of stability and non-stability.
The graph of the dependence of such boundary frequency on dimensionless parameter $\nu$ is shown
at Figure~\ref{fig:AM1:w-nu}, where the first fundamental frequency is also present
(the horizontal line).
The graph is constructed with the following dimensionless parameters:
\begin{equation}
\label{eq:param}
 \varepsilon_1 = 0.005, \quad \mu = 0.008, \quad \eta=7, \quad \delta=0.1 \, ,
\end{equation}
which take place in the real hybrid systems \cite{MR1835203}.

The intersection point in Figure~\ref{fig:AM1:w-nu} corresponds
to the value $\nu = 0.0529760481$, which is the maximum allowed value:
for all greater values of $\nu$, the natural frequency $\omega_1$
is in the non-stability domain that corresponds
to excitation of vibrations in the system.

As it can be easily seen from Figure~\ref{fig:AM1:w-nu}, at parameters \eqref{eq:param}
the vibrations will never take place on the second frequency $\omega_2$.

It is also possible to calculate the same boundary frequency, using the second asymptotical
method. To do this, we use the condition of equality in \eqref{eq:523}.
The results of such calculations at the same parameters \eqref{eq:param}
support the results in Figure~\ref{fig:AM1:w-nu}.

Now, let us also consider, at the same parameters \eqref{eq:param},
the results of calculations of the eigenvalues
of the nonconservative boundary value problem \eqref{eq:9}--\eqref{eq:11}
by the numerical method.

\begin{figure}
\subfigure[]{
  \label{fig:NM1:q1}
  \begin{minipage}[b]{0.5\textwidth}
  \centering
   \psfrag{nu}{$\nu$}
   \psfrag{q1}{$q_1$}
   \includegraphics[scale=0.4]{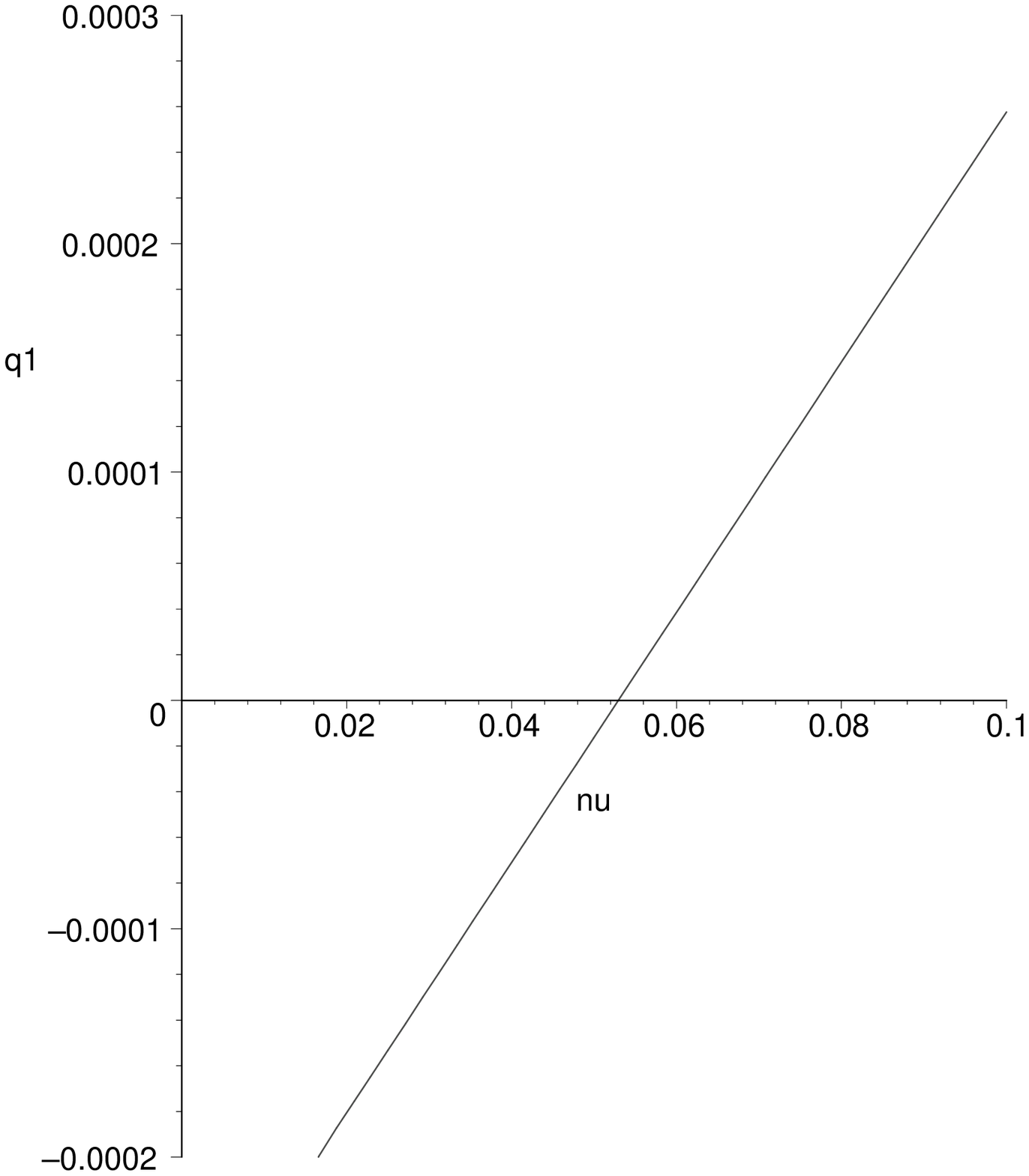}
  \end{minipage}}
\subfigure[]{
  \label{fig:NM2:q2}
  \begin{minipage}[b]{0.5\textwidth}
  \centering
 \psfrag{nu}{$\nu$}
 \psfrag{q2}{\hspace*{-0.2cm}$q_2$}
\includegraphics[scale=0.4]{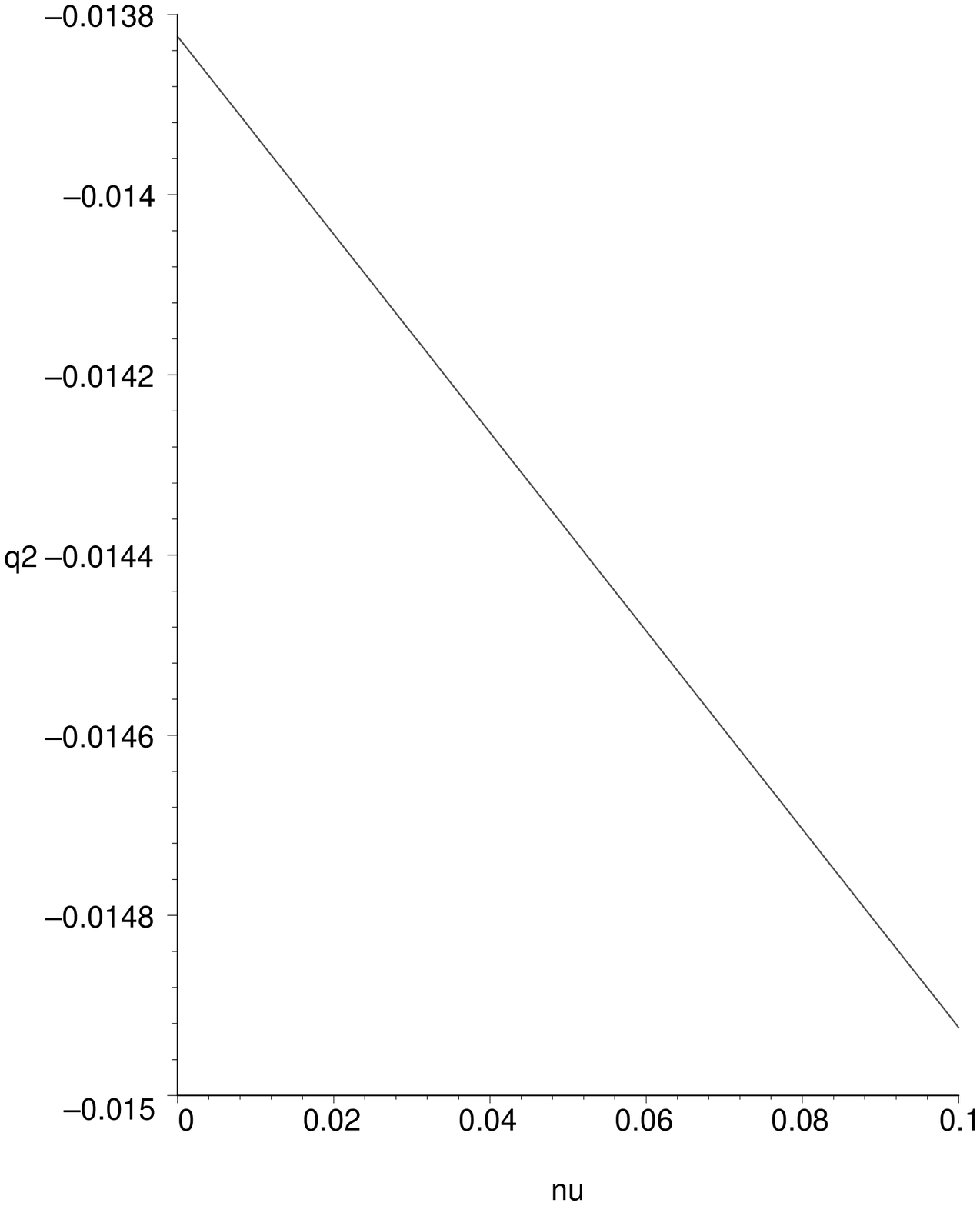}
  \end{minipage}}
\caption{Numerical Method: $\varepsilon_1 = 0.005$, $\mu = 0.008$, $\eta=7$, $\delta=0.1$}
\label{fig:NM}
\end{figure}

According with \eqref{eq:determ:div}, in Figure~\ref{fig:NM1:q1}
and Figure~\ref{fig:NM2:q2} we present, respectively,
the graphs of the dependence of real parts $q_{1}$ and $q_{2}$,
of the first two eigenvalues, on dimensionless parameter $\nu$.
According with \eqref{eq:determ:div},
the imaginary parts $\omega_1$, $\omega_2$ of the first two eigenvalues
are practically independent on parameter $\nu$, and coincide
with the results \eqref{eq:w}, obtained by the asymptotical method.

It is clear that the stability domain of the system corresponds to the negative values
of real parts $q_{1}$ and $q_{2}$ of the eigenvalues, and when they vanish,
we have the case of self-excited vibrations. So, from Figure~\ref{fig:NM2:q2}, we can see
that vibrations will not be excited at the second frequency. Besides, from Figure~\ref{fig:NM1:q1},
it is possible to calculate the value of parameter $\nu$ at which $q_{1}=0$. This value
of $\nu$ coincides, with accuracy $10^{-3}$, with the value obtained from the intersection
point in Figure~\ref{fig:AM1:w-nu} by asymptotical method.

Thus, all the results of calculations by the asymptotical and numerical methods
coincide with a sufficient accuracy.

In this way, for each definite hybrid system it is possible to define a character of possible vibrations.
As dimensionless parameter $\nu$ is directly proportional to the feedback coefficient $d$,
which may be purposefully changed by application of different feedbacks \cite{MR1835203}, it is possible to change
purposefully eigenvalues of the considered nonconservative boundary problem.
This allows to avoid in the system an excitation of undesirable vibrations, including
self-excited vibrations.

\section*{Acknowledgments}

The first author was supported by FCT (\emph{The Portuguese Foundation for Science
and Technology}), fellowship SFRH/BPD/14946/2004.
The authors are also grateful for the support of the R\&D unit
\emph{Centre for Research in Optimization and Control} (CEOC).



\bigskip
\bigskip


\begin{wrapfigure}[8]{l}{2.4cm}
\includegraphics[scale=0.5]{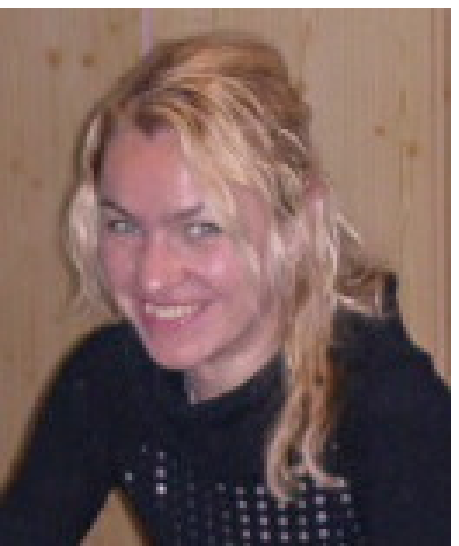}
\end{wrapfigure}
\noindent \textbf{Olena V. Mul} was born in 1972 in Ternopil, Ukraine.
She graduated in Applied Mathematics
from State University Lvivs'ka Politekhnica in Lviv,
Ukraine, in 1994. From 1994 to 1998 she was a post-graduate
student in the Institute of Cybernetics of Kyiv. She obtained her
PhD in Physics and Mathematics in 2001 from the Institute of Space
Research in Kyiv, Ukraine. She is with Ternopil State Technical
University, Ternopil, Ukraine, since 1995: as Engineer in 1995,
Assistant in 1996--2001, Associate Professor in 2001--2004. Since
2000 she is also with the Ternopil Academy of National Economy,
Ternopil, Ukraine: as an Assistant in 2000, as Senior Lecturer in
2001, and as an Associate Professor since 2002. She served as
Editor-in-Chief Deputy of International Scientific Journal
\emph{Computing} in 2002--2004. In 2004 she got a Research
Fellowship in the Department of Mathematics of the University of
Aveiro, Portugal, from the Portuguese Foundation for Science and
Technology (FCT).

\bigskip
\bigskip

\begin{wrapfigure}[8]{l}{2.4cm}
\includegraphics[scale=0.15]{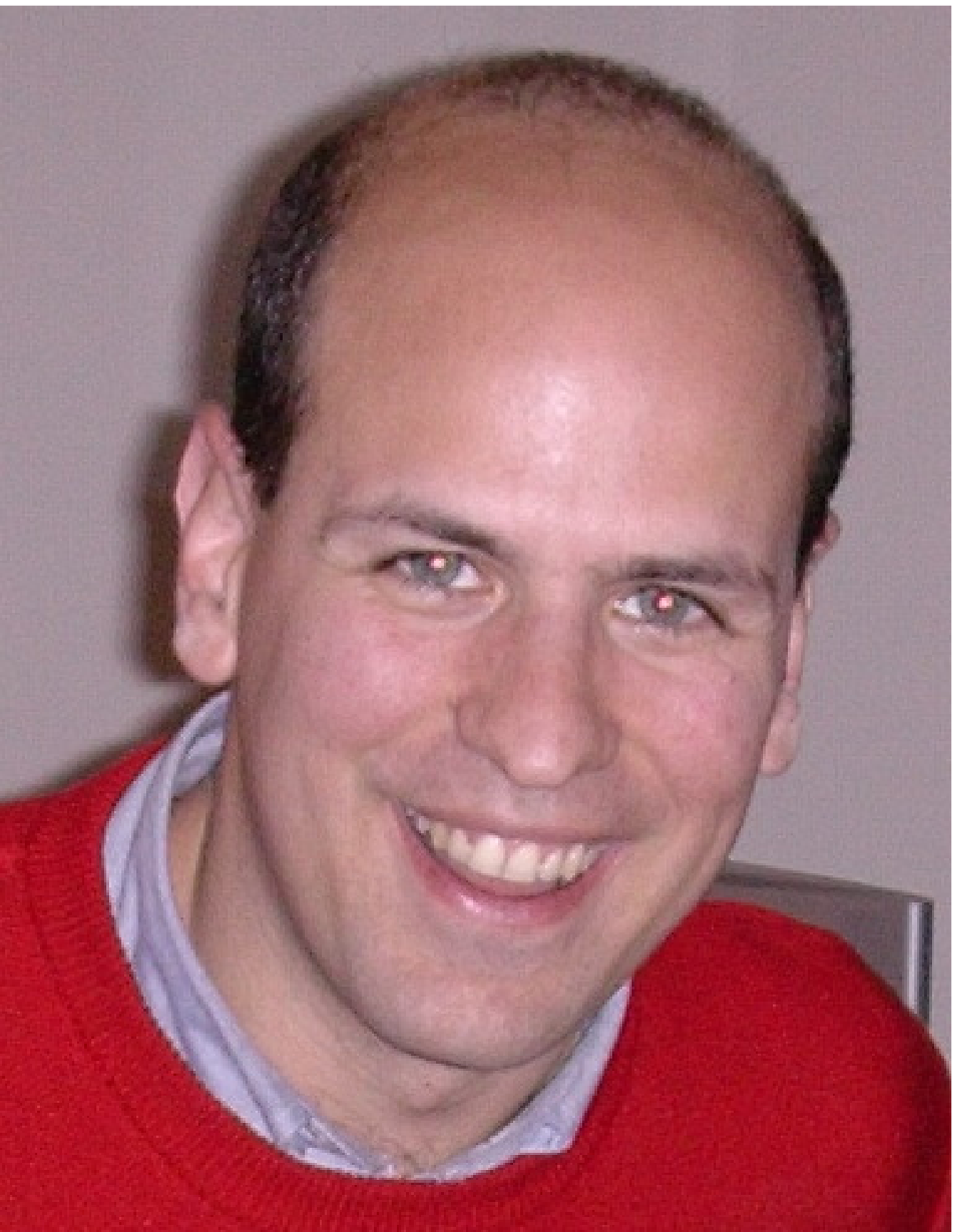}
\end{wrapfigure}
\noindent \textbf{Delfim F. M. Torres}
was born in 1971 in Mozambique, Africa. He graduated from
University of Coimbra, Portugal, in Computer Science Engineering
in 1994. From 1993 to 1994 he was a Junior Researcher at the
Institute of Biomedical Research on Light and Image, in Coimbra,
Portugal. In 1994--1998 he did his master studies at the
Department of Mathematics of the University of Aveiro, Portugal,
where he held a Teaching Assistant position. February 1998 he
obtained his MSc in Mathematics, specialization on Optimization
and Control Theory, and got a position of Assistant Lecturer. He
obtained his PhD in Mathematics in 2002 from the University of
Aveiro, Portugal, where he is now an Assistant Professor. He
served as vice-president of the Department of Mathematics of the
University of Aveiro in 2002--2004. Since 2004 he is the
Scientific Coordinator of the Control Theory Group (\textsf{cotg})
of the Centre for Research in Optimization and Control
(\textsf{CEOC}) of the University of Aveiro. His scientific CV can
be found at \texttt{<http://www.mat.ua.pt/delfim>}.



\begin{bibdiv}

\begin{biblist}

\bib{MR1835203}{article}{
    author={Krasnoshapka, V. A.},
     title={On the stability of systems with distributed parameters in the
            realization of the active control of large flexible systems},
  language={Russian},
   journal={Problemy Upravlen. Inform.},
      date={2000},
    number={6},
     pages={60\ndash 65, 154},
      issn={0572-2691},
    review={MR 1835203},
}

\bib{MR1723682}{article}{
    author={Krasnoshapka, V. A.},
    author={Mul, O. V.},
     title={On the excitation of auto-oscillations in controlled dynamical
            systems with distributed parameters},
  language={Russian},
   journal={Problemy Upravlen. Inform.},
      date={1999},
    number={3},
     pages={73\ndash 79, 153},
      issn={0572-2691},
    review={MR 1723682},
}

\bib{Norm-fund}{book}{
    author={Kuhta, K. Ya.},
    author={Kravchenko, V. P.},
     title={Normalnye fundamentalnye sistemy v zadachakh teorii kolebanij
            [Normal Fundamental Systems in Vibration Theory Problems]},
  language={Russian},
 publisher={``Naukova Dumka''},
     place={Kiev},
      date={1973},
     pages={205},
}

\bib{MR80f:73002}{book}{
    author={Kuhta, K. Ya.},
    author={Kravchenko, V. P.},
     title={Nestatsionarnye granichnye zadachi s nepreryvno-diskretnymi parametrami
            [Nonstationary boundary value problems with continuous-discrete parameters]},
  language={Russian},
 publisher={``Naukova Dumka''},
     place={Kiev},
      date={1978},
     pages={217},
    review={MR 80f:73002},
}

\bib{MR88k:49009}{book}{
    author={Kuhta, K. Ya.},
    author={Kravchenko, V. P.},
    author={Krasnoshapka, V. A.},
     title={Kachestvennaya teoriya upravlyaemykh dinamicheskikh sistem s
            nepreryvno-diskretnymi parametrami
             [The qualitative theory of controllable dynamical systems with continuous-discrete parameters]},
  language={Russian},
 publisher={``Naukova Dumka''},
     place={Kiev},
      date={1986},
     pages={223},
    review={MR 88k:49009},
}

\bib{MR1383822}{book}{
    author={Landa, P. S.},
     title={Nonlinear oscillations and waves in dynamical systems},
    series={Mathematics and its Applications},
    volume={360},
 publisher={Kluwer Academic Publishers Group},
     place={Dordrecht},
      date={1996},
     pages={xvi+538},
      isbn={0-7923-3931-2},
    review={MR1383822 (97h:00008)},
}

\bib{myNY98}{article}{
    author={Mul, O.},
     title={On Conditions of Excitation of Self-Oscillations
            in a Nonconservative Dynamic System with Distributed Parameters},
   journal={Cybernetics and Computing Technology, Complex Control Systems,
            Allerton Press, New York},
      date={1998},
    number={111},
     pages={70\ndash72},
}

\bib{myDresd}{article}{
    author={Mul, O.},
    author={Kravchenko, V.}
     title={Investigations of Vibrations in the
            Complex Dynamical Systems of Transmission Pipelines,
            "Interface and Transport Dynamics. Computational Modelling"},
   journal={Lecture Notes in Computational Science and Engineering,
            Springer-Verlag Berlin Heidelberg},
      date={2003},
    number={32},
     pages={295\ndash300},
}

\bib{MR2091174}{book}{
    author={Samoilenko, A.},
    author={Petryshyn, R.},
     title={Multifrequency oscillations of nonlinear systems},
    series={Mathematics and its Applications},
    volume={567},
 publisher={Kluwer Academic Publishers Group},
     place={Dordrecht},
      date={2004},
     pages={vi+317},
      isbn={1-4020-2030-9},
    review={MR2091174},
}

\bib{MR2081767}{book}{
    author={Svetlitsky, V. A.},
     title={Engineering vibration analysis: worked problems. 1},
    series={Foundations of Engineering Mechanics},
      note={Translated from the Russian by G. I. Merzon and V. A. Chechin},
 publisher={Springer-Verlag},
     place={Berlin},
      date={2004},
     pages={x+316},
      isbn={3-540-20658-2},
    review={MR2081767},
}

\end{biblist}

\end{bibdiv}
\end{document}